\documentclass[letterpaper, paper,11pt]{AAS}		% for final proceedings (20-page limit)

\usepackage{bm}
\usepackage{amsmath,mathtools,upgreek,gensymb}
\usepackage{subfigure}
\usepackage[colorlinks=true, pdfstartview=FitV, linkcolor=black, citecolor= black, urlcolor= black]{hyperref}
\usepackage{overcite}
\usepackage{footnpag}			      	% make footnote symbols restart on each page

\PaperNumber{19-810}

\begin{document}

\title{Second-Order Solution for Relative Motion on Eccentric Orbits in Curvilinear Coordinates}

\author{Matthew Willis\thanks{PhD Candidate, Department of Mechanical Engineering, Stanford University, 496 Lomita Mall, Stanford, CA 94305.}, Kyle T. Alfriend\thanks{Professor, Department of Aerospace Engineering, Texas A\&M University, 3141 TAMU, College Station, TX 77843}, and Simone D'Amico\thanks{Professor, Department of Aeronautics and Astronautics, Stanford University, 496 Lomita Mall, Stanford, CA 94305.}
}

\maketitle{}

\begin{abstract}
A new, second-order solution in curvilinear coordinates is introduced for the relative motion of two spacecraft on eccentric orbits. The second-order equations for unperturbed orbits are derived in spherical coordinates with true anomaly as the independent variable, and solved by the method of successive approximations. A comparison of error trends against eccentricity and inter-spacecraft separation is presented between the new solution and prominent Cartesian, curvilinear, and orbital element based solutions from the literature. The second-order curvilinear solution offers a thousand-fold improvement in accuracy over the first-order curvilinear solution, and still greater improvement over first- and second-order rectilinear solutions when large along-track separations are present.

. 
\end{abstract}

\section{Introduction}

Distributed space systems are a mission-enabling technology for commercial and scientific applications ranging from on-orbit satellite inspection and servicing to observations of gravitational waves and direct imaging of extrasolar planets.\cite{DistributedSS} Advanced formation guidance, navigation, and control algorithms are needed to make such missions a reality, and these will rely heavily on the model used to describe the dynamics of spacecraft relative motion.
The limited processing power typical of flight hardware introduces a tradeoff between computational efficiency and model accuracy for onboard implementation. Analytical solutions are particularly valuable because their accuracy is not tied to an integration step size or iteration tolerance, and therefore does not scale uniformly with computational cost. There are two broad categories of relative motion models: those based on orbital elements and those using a translational state representation. Sullivan and D'Amico conducted a thorough survey of existing dynamics models and solutions in both categories, including a comparison of their performance under various assumptions.\cite{SullivanGrimberg} Orbital element representations offer better accuracy due to their fundamental connection to the underlying physics and relative motion geometry. However, spacecraft sensors and actuators do not live in orbital element space, so there is an advantage to using translational state models that avoid this intermediate representation. The present work is therefore focused on high-fidelity solutions in this category.

\begin{figure}[b!]
	\centering\includegraphics[width=\textwidth,trim=0in 2in 2.2in 0, clip]{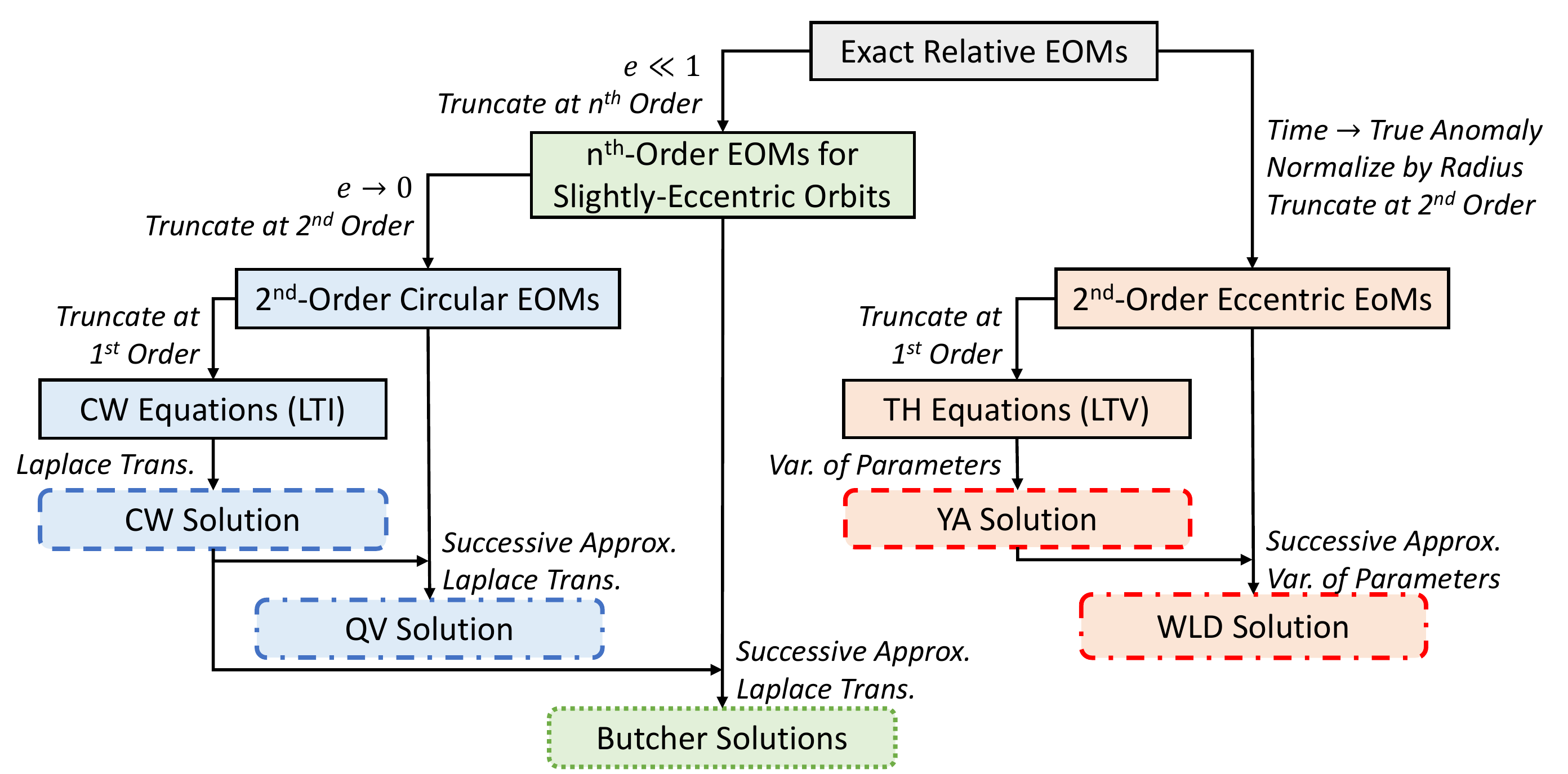}
	\caption{Derivation family tree of translational state solutions.}
	\label{fig:methods}
\end{figure}

A family tree of translational state models and solutions is presented in Figure~\ref{fig:methods}, emphasizing how each is obtained from the equations of motion (EOMs).
The best-known of these solutions is that of Clohessy and Wiltshire (CW), which addresses the linear, time-invariant problem of relative motion between two spacecraft in close proximity on near-circular orbits.\cite{HCW} A second-order solution to the circular orbit problem was independently derived by London and Sasaki and later by Stringer and Newman, and is often referred to as the Quadratic-Volterra (QV) solution.\cite{London,SasakiQV,NewmanQV,QVcompare} These authors obtained the second-order solution using the method of successive approximations, wherein the first-order CW solution is substituted into the nonlinear dynamics, resulting in an inhomogeneous linear system that may be solved by elementary differential equations techniques. The same strategy has been used by Melton and Butcher, et al. to incorporate leading-order effects of eccentricity.\cite{Melton, ButcherLovellCart} 
%Because these solutions are explicit in time, they avoid the need to numerically solve Kepler's equation for true anomaly at fixed times.\cite{Vallado}
Tschauner and Hempel (TH) provided a linear description of the dynamics governing relative motion on elliptical orbits by appropriately normalizing the coordinates and changing the independent variable from time to true anomaly.\cite{TH} Solutions to this system by Tschauner and Hempel, Carter, and others offered better accuracy than the CW solution in slightly eccentric orbits but suffered from singularities at zero eccentricity.\cite{CarterTH} Yamanaka and Ankersen (YA) were able to remove this singularity with the use of a new integral that grows in proportion to time.\cite{YA} Willis, Lovell, and D'Amico (WLD) recently introduced a second-order solution for relative motion on eccentric orbits by applying the method of successive approximations to the first-order YA solution.\cite{WillisLovellDAmico} 
The present work extends this contribution by introducing a second-order solution to the eccentric orbit problem in curvilinear coordinates. While this paper will focus on the closely-related models and solutions shown in Figure~\ref{fig:methods}, other approaches exist to arrive at similar translational state solutions. For example, the higher-order state transition tensor theory studied by Park and Scheeres could be used in place of differential equations techniques to develop an equivalent second-order relative motion solution.\cite{ScheeresSTT}

% Talk about curvilinear solutions
The first-order equations of relative motion for circular orbits have the same form in curvilinear coordinates as in rectilinear coordinates, and are therefore solved by CW.\cite{AlfriendFF} For large along-track separations, the rectilinear model gives a poor approximation of the relative state whereas the curvilinear model naturally captures the orbit geometry and enables more accurate propagation. Butcher et al. exploited the curvilinear version of CW to develop higher-order solutions in spherical and cylindrical coordinates for circular and slightly-elliptical orbits.\cite{ButcherBurnettCurv} While the first-order dynamics are identical in rectilinear and curvilinear coordinates, the higher-order dynamics are fundamentally different. Thus, a second-order solution in spherical coordinates is not the same as a second-order Cartesian solution after a coordinate transformation.\cite{ButcherLovellDT} It has long been suspected that like the CW solution, there is a spherical coordinate equivalent of YA with identical form. Han et al. recently demonstrated that the TH equations are identical in Cartesian and spherical coordinates and that the YA solution is therefore valid in both, laying the foundation for the second-order solution derived herein.\cite{HanYASpherical}

%% NEED TO REVISE
The body of this paper is divided into three sections. First, the relevant background material is reviewed, including the spherical coordinate definitions and their relation to the more familiar Cartesian coordinates, the development of the second-order equations of relative motion, and the introduction of the YA solution. This is followed by the derivation of the new, second-order solution in curvilinear coordinates. The derivation details the method of successive approximations and concludes with the explicit statement of the new solution. In the third section, the solution is validated through a performance comparison with a selection of rectilinear and curvilinear solutions from the literature. The paper concludes with a brief summary of the results and discussion of future directions for research.

\section{Background}

\subsection{Cartesian vs. Spherical Coordinates}

\begin{figure}[bth]
    \centering
    \includegraphics[width=2.5in,trim=0.7in 0in 0.3in 0, clip]{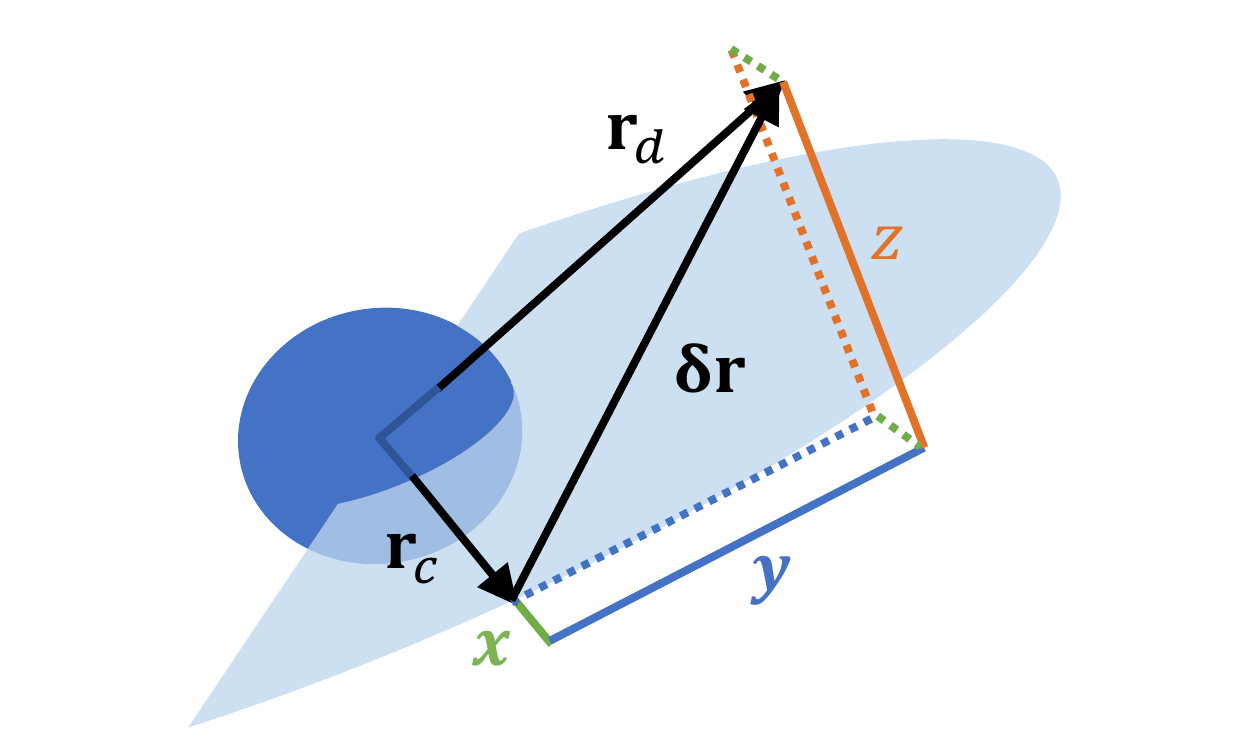}
    \includegraphics[width=2.5in,trim=0.3in 0in 0.7in 0, clip]{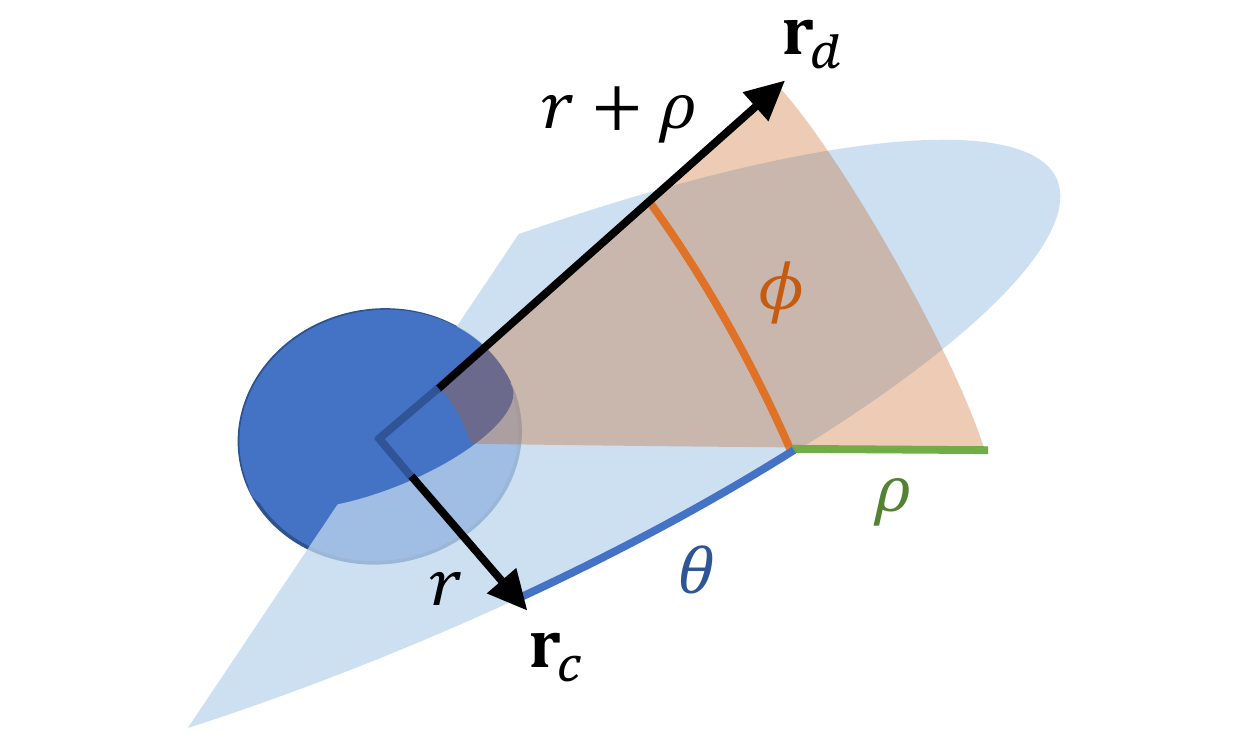}
    \caption{Definition of Cartesian (left) and spherical coordinates (right) for relative motion.}
    \label{fig:coords}
\end{figure}

Cartesian coordinates provide a convenient means of describing the position vector $\mathbf{\boldsymbol \updelta r}$ of a deputy spacecraft relative to a chief located at $\mathbf{r}_c$, a distance $r$ from the central body. We express the relative motion with respect to the Radial-Transverse-Normal (RTN) frame rotating with the chief's orbit using the associated $x$, $y$, and $z$ coordinates. The $x$ axis extends radially away from the central body, $z$ extends along the direction normal to the chief's orbital plane, and $y$ completes the orthogonal basis with positive component in the direction of motion. As illustrated on the left of Figure~\ref{fig:coords}, the relative position vector is $\boldsymbol{\updelta}\mathbf{r} = [x, y, z]^T$.

Many curvilinear coordinate systems may be used to describe the relative motion of two spacecraft, but this paper is concerned only with the spherical coordinates $(\rho,\theta,\phi)$, illustrated on the right of Figure~\ref{fig:coords}. These are defined so that $\rho = r_d - r$ is the difference in radial separation from the central body between the two spacecraft, $\theta$ is the angle from the chief's position vector to the projection of the deputy's position vector onto the chief's orbital plane and $\phi$ is the angle from this projection to the deputy's position vector.

Figure~\ref{fig:coords} emphasizes the analogous roles of the Cartesian coordinates $(x,y,z)$ and spherical coordinates $(\rho,\theta,\phi)$ for describing radial, along-track, and cross-track separations, respectively. For small separations, the Cartesian and curvilinear coordinates are related by $\boldsymbol{\updelta}\mathbf{r} = [x, y, z]^T \approx [\rho, r\theta, r\phi]^T$. However, it is also apparent in the figure that this relationship breaks down for larger separations. The exact mapping from spherical to Cartesian coordinates is given by
\begin{equation}\label{eq:curv2cart}
    \begin{aligned}
        x &= (r + \rho) \cos \phi \cos \theta - r \\
        y &= (r + \rho) \cos \phi \sin \theta \\
        z &= (r + \rho) \sin \phi \\
    \end{aligned}
\end{equation}
and a complete set of transformations between the spherical coordinates and the relative position and velocity vectors are provided in the appendix. Although the two descriptions are closely related, it is important to recognize that the Cartesian coordinates directly describe the evolution of the relative position and velocity \textit{vectors} with respect to the rotating RTN frame, while the spherical coordinates use time-derivatives of \textit{scalar} quantities to describe the difference in inertial motion of the chief and deputy. 

To derive the equations of relative motion in spherical coordinates, we begin with the acceleration of each spacecraft in an inertial frame. The chief's motion is governed by
\begin{equation}\label{eq:eom_chief_vec}
    \Ddot{r} \hat{\mathbf{r}}_c + 2 \dot{r} {}^I \mathbf{\boldsymbol \upomega}^c \times \hat{\mathbf{r}}_c + r \left( {}^I\dot{\mathbf{\boldsymbol \upomega}}^c \times \hat{\mathbf{r}}_c + {}^I \mathbf{\boldsymbol \upomega}^c \times \left( {}^I \mathbf{\boldsymbol \upomega}^c \times \hat{\mathbf{r}}_c \right) \right) = -\frac{\mu}{r^2} \hat{\mathbf{r}}_c
\end{equation}
where $\hat{\mathbf{r}}_c$ is the unit vector from the central body to the chief's position and ${}^I \mathbf{\boldsymbol \upomega}^c$ is the angular velocity of its RTN frame with respect to the inertial frame.
In the absence of perturbations, the vectors are expressed in RTN components as $\hat{\mathbf{r}}_c = [r, 0, 0]^T$ and ${}^I\dot{\mathbf{\boldsymbol \upomega}}^c = [0,0,\dot{\theta}_c]$, where $\theta_c$ is the angle from an arbitrary reference vector in the chief's orbital plane to the chief's position vector $\mathbf{r}_c$. Equation~(\ref{eq:eom_chief_vec}) can then be converted to the scalar relations
\begin{equation}\label{eq:eom_chief}
    \begin{aligned}
        \Ddot{r} - r\dot{\theta}_c^2 &= -\frac{\mu}{r^2} \\
        2\dot{r}\dot{\theta}_c + r \Ddot{\theta}_c &= 0
    \end{aligned}
\end{equation}
Similar to Equation~(\ref{eq:eom_chief_vec}), the evolution of the deputy's position vector is described by
%{}^I \Ddot{\mathbf{r}}_d = 
\begin{equation}\label{eq:eom_dep}
    (\Ddot{r}+\Ddot{\rho}) \hat{\mathbf{r}}_d + 2(\dot{r}+\dot{\rho}) {}^I \mathbf{\boldsymbol \upomega}^d \times \hat{\mathbf{r}}_d + (r + \rho) \left( {}^I\dot{\mathbf{\boldsymbol \upomega}}^d \times \hat{\mathbf{r}}_d + {}^I \mathbf{\boldsymbol \upomega}^d \times \left( {}^I \mathbf{\boldsymbol \upomega}^d \times \hat{\mathbf{r}}_d \right) \right) = -\frac{\mu}{(r + \rho)^2} \hat{\mathbf{r}}_d
\end{equation}
We relate the deputy's direction vector $\hat{\mathbf{r}}_d$ and orbital angular velocity ${}^I \mathbf{\boldsymbol \upomega}^d$ to our curvilinear coordinates by expressing Equation~(\ref{eq:eom_dep}) in the chief's RTN components. The vectors are given by (cf. Equation~\ref{eq:curv2cart})
\begin{equation}\label{eq:dep_vecs}
    \hat{\mathbf{r}}_d = \begin{bmatrix}
        \cos \phi \cos \theta \\
        \cos \phi \sin \theta \\
        \sin \phi
    \end{bmatrix} \quad {}^I \mathbf{\boldsymbol \upomega}^d = \begin{bmatrix}
        \dot{\phi} \sin \theta \\
        -\dot{\phi} \cos \theta \\
        \dot{\theta} + \dot{\theta}_c
    \end{bmatrix} \quad {}^I \dot{\mathbf{\boldsymbol \upomega}}^d = \begin{bmatrix}
        \Ddot{\phi} \sin \theta + (\dot{\theta}+\dot{\theta}_c) \dot{\phi} \cos \theta \\
        -\Ddot{\phi} \cos \theta + (\dot{\theta} + \dot{\theta}_c) \dot{\phi} \sin \theta \\
        \Ddot{\theta} + \Ddot{\theta}_c
    \end{bmatrix}
\end{equation}
Substituting Equation~(\ref{eq:dep_vecs}) into Equation~(\ref{eq:eom_dep}) and solving for $\Ddot{\rho}$, $\Ddot{\theta}$, and $\Ddot{\phi}$ leads to the system
\begin{equation}\label{eq:eom_curv_exact}
    \begin{aligned}
        \Ddot{\rho} &= -\Ddot{r} - \frac{\mu}{(r + \rho)^2} + (r + \rho)\left(\dot{\phi}^2 + (\dot{\theta} + \dot{\theta}_c)^2\cos^2 \phi \right) \\
        \Ddot{\theta} &= -\Ddot{\theta}_c + 2(\dot{\theta} + \dot{\theta}_c) \dot{\phi} \tan \phi - 2\frac{(\dot{r} + \dot{\rho})}{(r+\rho)}(\dot{\theta} + \dot{\theta}_c) \\
        \Ddot{\phi} &=  -2\frac{(\dot{r}+\dot{\rho})}{(r+\rho)} \dot{\phi} - (\dot{\theta} + \dot{\theta}_c)^2 \cos \phi \sin \phi
    \end{aligned}
\end{equation}
Up to this point, no approximations have been made beyond the assumption of unperturbed Keplerian motion. We wish to approximate the equations of motion as a system of polynomials in the relative state variables. Series expanding the nonlinearities and truncating at second order, Equation~(\ref{eq:eom_curv_exact}) becomes
\begin{equation}\label{eq:eom2_curv_dt}
    \begin{aligned}
        \Ddot{\rho} &= 2\frac{\mu}{r^3} \rho - 3 \frac{\mu}{r^4} \rho^2 + \rho \dot{\theta}_c^2 + 2 r \dot{\theta}_c \dot{\theta} + r \dot{\phi}^2 + 2  \dot{\theta}_c \rho \dot{\theta} + r \dot{\theta}^2 - r \dot{\theta}_c^2 \phi^2 \\
        \Ddot{\theta} &= -\frac{1}{r}\left(\rho \Ddot{\theta}_c + 2 \dot{\rho} \dot{\theta}_c + 2 \dot{r}\dot{\theta} + 2 \dot{\rho} \dot{\theta} - 2 r \dot{\theta}_c \dot{\phi}\phi - \frac{\Ddot{\theta}_c}{r} \rho^2 - 2\frac{\dot{\theta}_c}{r} \rho \dot{\rho} - 2\frac{ \dot{r}}{r} \rho \dot{\theta} \right) \\
        \Ddot{\phi} &= - \frac{1}{r}\left(r \dot{\theta}_c^2 \phi + 2 \dot{r} \dot{\phi} + 2 r \dot{\theta}_c \dot{\theta} \dot{\phi} + 2 \dot{\rho} \dot{\phi} - 2\frac{\dot{r}}{r} \rho \dot{\phi} \right)
    \end{aligned}
\end{equation}

Following the approach of Tschauner and Hempel, we nondimensionalize the equations by changing the independent variable from time to true anomaly $f$ and normalizing the coordinates by the chief's orbit radius $r$. For convenience, we introduce the parameter
\begin{equation}
    k = \frac{p}{r} = 1 + e \cos{f}
\end{equation}
where $p = a(1-e^2)$ is the semi-latus rectum. In previous work by the authors and in the work of Yamanaka and Ankersen the parameter $k$ was denoted by $\rho$. This work adopts the notation of Alfriend et al. to avoid confusion with the curvilinear coordinates.\cite{AlfriendFF} 
We will denote the normalized coordinates with $(\Tilde{\;})$ and derivatives with respect to true anomaly by $(\;)'$. 
Because the angular coordinates $\theta$ and $\phi$ are already nondimensional, only $\rho$ must be normalized according to $\Tilde{\rho} = \rho / r$. The independent variable is changed to true anomaly using the substitutions
\begin{equation}
    \begin{aligned}
        \dot{r} &= r\dot{\theta}_c \frac{e \sin f}{k} \\
        \dot{\rho} &= \frac{r \dot{\theta}_c}{k} \left(k \Tilde{\rho}' + \Tilde{\rho}e \sin f \right) \\
        \Ddot{\rho} &= \frac{r \dot{\theta}_c^2}{k} \left(k \Tilde{\rho}'' + \Tilde{\rho}' e \cos f\right) \\
        \dot{\alpha} &= \alpha' \dot{\theta}_c \\
        \Ddot{\alpha} &= \left( \alpha'' - 2 \alpha' \frac{e \sin f}{k} \right)\dot{\theta}_c^2
    \end{aligned}
\end{equation}
where $\alpha$ is either $\theta$ or $\phi$. Equation~(\ref{eq:eom_chief}) can then be used to eliminate $r$, $\theta_c$, and their derivatives. Performing these transformations on Equation~(\ref{eq:eom2_curv_dt}) leads to the second-order system we wish to solve,
\begin{equation}\label{eq:eom2_curv}
    \begin{aligned}
        \Tilde{\rho}'' - 2\theta' - \frac{3}{k} \Tilde{\rho} &= -\frac{3}{k} \Tilde{\rho}^2 + 2 \Tilde{\rho} \theta' + \phi'^2 + \theta'^2 - \phi^2 \\
        \theta'' + 2\Tilde{\rho}' &= -2 \Tilde{\rho}' \theta' + 2 \phi' \phi + 2 \Tilde{\rho} \Tilde{\rho}' \\
        \phi'' + \phi &= -2 \theta' \phi - 2 \Tilde{\rho}' \phi'
    \end{aligned}
\end{equation}
In Equation~(\ref{eq:eom2_curv}), the first-order terms that appear in the Tschauner-Hempel equations have been moved to the left while the nonlinear, second-order terms remain on the right-hand side. 

Before proceeding to the first- and second-order solutions, it is worth recalling the equations of relative motion in Cartesian coordinates for comparison. A detailed development of the second-order Cartesian equations is provided in the original derivation of the WLD solution, but here we will reproduce only the fundamental description of the relative dynamics and final second-order system.\cite{WillisLovellDAmico} Let ${}^c (\dot{\;})$ denote a time derivative with respect to the chief's RTN frame, $\boldsymbol{\updelta}\mathbf{v} \equiv {}^c \boldsymbol{\updelta}\dot{\mathbf{r}}$ be the relative velocity, and ${}^I \boldsymbol{\upomega}^c$ be the angular velocity of the rotating frame with respect to the inertial frame. After applying the theorem of Coriolis to the fundamental orbital differential equations of chief and deputy and expanding the differential gravitational attraction to second-order in $(\delta r/r)$, the equations of relative motion in the chief's RTN frame are given by
\begin{small}
\begin{equation}\label{eq:eom_full}
    {}^c \boldsymbol{\updelta} \Ddot{\mathbf{r}} = -\frac{\mu}{r^3} \left(\mathbf{ \boldsymbol{\updelta}r} - 3 \frac{\mathbf{r \boldsymbol{\cdot} \boldsymbol{\updelta} r}}{r^2} (\mathbf{r + \boldsymbol{\updelta}r}) - \frac{3}{2} \frac{\delta r^2}{r^2} \mathbf{r} + \frac{15}{2} \frac{(\mathbf{r \boldsymbol{\cdot} \boldsymbol{\updelta} r})^2}{r^4} \mathbf{r} \right) - 2 {}^I\boldsymbol{\upomega}^c \times \boldsymbol{\updelta} \mathbf{v} - {}^I\Dot{\boldsymbol{\upomega}}^c \times \boldsymbol{\updelta} \mathbf{r} - {}^I\boldsymbol{\upomega}^c \times {}^I\boldsymbol{\upomega}^c \times \boldsymbol{\updelta} \mathbf{r}
\end{equation}
\end{small}

\vspace{-1em}
After series-expanding the nonlinearities and nondimensionalizing the system, the second-order equations of motion take on the form
\begin{equation}\label{eq:eom2_cart}
    \begin{aligned}
        \Tilde{x}'' - 2\Tilde{y}' - \frac{3}{k} \Tilde{x} &= -\frac{3}{k} \Tilde{x}^2 + \frac{3}{2 k} (\Tilde{y}^2 + \Tilde{z}^2 ) \\
        \Tilde{y}'' + 2\Tilde{x}' &= \frac{3}{k} \Tilde{x} \Tilde{y} \\
        \Tilde{z}'' + \Tilde{z} &= \frac{3}{k} \Tilde{x} \Tilde{z}
    \end{aligned}
\end{equation}
As in Equation~(\ref{eq:eom2_curv}), the first-order terms have been collected on the left-hand side of Equation~(\ref{eq:eom2_cart}).

Comparing the curvilinear system in Equation~(\ref{eq:eom2_curv}) with its Cartesian counterpart in Equation~(\ref{eq:eom2_cart}) reveals several remarkable similarities and differences. First, as demonstrated by Han et al., the linear terms are identical in form and accept the same solutions.\cite{HanYASpherical} The Yamanaka-Ankersen state transition matrix introduced below can therefore be used for either system with the substitutions $\Tilde{x} \leftrightarrow \Tilde{\rho}$, $\Tilde{y} \leftrightarrow \theta$, and $\Tilde{z} \leftrightarrow \phi$. However, the second-order terms capture different nonlinearities in the two coordinate systems. The first-order terms in the rectilinear coordinates accurately describe the fictitious forces due to the rotating reference frame---the terms appearing outside the parentheses in Equation~(\ref{eq:eom_full}). This stems from the fact that the axis of rotation is fixed in the RTN coordinate system. As a result, the second-order terms all arise from approximation of the differential gravitational effect on the two spacecraft. The difference in radial separation is seen in the $\Tilde{x}$ equation and corrections for the difference in the radial directions are seen in all three equations. In contrast, the curvilinear coordinates condense the differential gravity correction to a single term in the $\Tilde{\rho}$ equation. The other terms appearing in this formulation are due to the difference in chief and deputy angular velocity vectors, which govern the evolution of their respective orbit radii and direction angles about the central body. Technically, these are kinematic expressions stemming from the spherical coordinate description of the motion, but are closely related to the fictitious forces of a rotating reference frame. The non-gravitational terms in the $\Tilde{\rho}$ equation are linked to the centrifugal force, as is the $\theta' \phi$ term in the $\phi$ equation. All other terms are connected to the Coriolis force except $\phi' \phi$ in the $\theta$ equation, which combines effects from centrifugal and Euler forces.

\subsection{Yamanaka-Ankersen Solution}

Yamanaka and Ankersen found an analytical solution to the linear Tschauner-Hempel equations, obtained by dropping the right-hand side of Equation~(\ref{eq:eom2_cart}). An identical system in spherical coordinates results from dropping the right-hand side of Equation~(\ref{eq:eom2_curv}),
\begin{equation}\label{eq:TH}
    \begin{aligned}
        \Tilde{\rho}'' - 2\theta' - \frac{3}{k} \Tilde{\rho} &= 0 \\
        \theta'' + 2\Tilde{\rho}' &= 0 \\
        \phi'' + \phi &= 0
    \end{aligned}
\end{equation}
Their key contribution was to eliminate singularities in the solution to the TH equations by introducing the integral $J(t)$, defined as
\begin{equation}
    J(t) = \int_{f_0}^f \frac{d\tau}{k(\tau)^2} =  \sqrt{\frac{\mu}{p^3}} (t-t_0) % = \frac{M-M_0}{(1-e^2)^{3/2}}
\end{equation}
Although the integration is taken over true anomaly, $J(t)$ is a linear function of time.
%\begin{align*}
%    \Tilde{x}_1 &= K_1 \left( 1 - \frac{3}{2} e k J(t) \sin{f} \right) + K_2 k \sin{f} + K_3 k \cos{f} \\
%    \Tilde{y}_1 &= K_4 + K_2 (1 + k) \cos{f} - K_3 (1 + k) \sin{f} - \frac{3}{2} K_1 k^2 J(t) \\
%    \Tilde{z}_1 &= K_5 \sin{f} + K_6 \cos{f}
%\end{align*}
The solution to Equation~(\ref{eq:TH}) is given by the linear system
\begin{equation}\label{eq:YA}
    \begin{bmatrix}
        \Tilde{\rho} \\
        \theta \\
        \phi \\
        \Tilde{\rho}' \\
        \theta' \\
        \phi'
    \end{bmatrix} = \begin{bmatrix}
        \left(1 - \frac{3}{2} e k J(t) \sin{f} \right) & k \sin{f} & k \cos{f} & 0 & 0 & 0 \\
        -\frac{3}{2} k^2 J(t) & (1 + k) \cos{f} & -(1 + k) \sin{f} & 1 & 0 & 0 \\
        0 & 0 & 0 & 0 & \sin{f} & \cos{f} \\
        - \frac{3}{2} e \left( (k \sin{f})' J(t) + \frac{\sin{f}}{k}\right) & (k \sin{f})' & (k \cos{f})' & 0 & 0 & 0 \\
        \frac{3}{2}(2 e k J(t) \sin{f} - 1) & -2 k \sin{f} & e - 2 k \cos{f} & 0 & 0 & 0 \\
        0 & 0 & 0 & 0 & \cos{f} & -\sin{f}
    \end{bmatrix} \begin{bmatrix}
        K_1 \\
        K_2 \\
        K_3 \\
        K_4 \\
        K_5 \\
        K_6
    \end{bmatrix}
\end{equation}
where $(k \sin{f})' = \cos{f} + e \cos{2f}$ and $(k \cos{f})' = -(\sin{f} + e \sin{2f})$, and $K_1$ through $K_6$ are integration constants. The relative velocity components are computed by differentiation of the corresponding relative position solution component with respect to true anomaly and do not contribute any additional solution information.

To express the solution in terms of initial conditions, one may solve for the integration constants by inverting Equation~(\ref{eq:YA}) and evaluating at the initial time $t_0$. Using $J(t_0) = 0$, this leads to
\begin{equation}\label{eq:YA_inv}
    \begin{bmatrix}
        K_1 \\
        K_2 \\
        K_3 \\
        K_4 \\
        K_5 \\
        K_6
    \end{bmatrix} = \begin{bmatrix}
        \frac{6 k_0 + 2 e^2 - 2}{1-e^2} & 0 & 0 & \frac{2 e k_0 \sin{f_0}}{1-e^2} & \frac{2 k_0^2}{1-e^2} & 0 \\
        -3 \left(1 + \frac{e^2}{k_0} \right) \frac{\sin{f_0}}{1-e^2} & 0 & 0 & \frac{k_0 \cos{f_0} - 2 e}{1-e^2} & -\frac{1 + k_0}{1-e^2} \sin{f_0} & 0 \\
        -3 \frac{e + \cos{f_0}}{1-e^2} & 0 & 0 & -\frac{k_0 \sin{f_0}}{1-e^2} & -\frac{e + (1 + k_0) \cos{f_0}}{1-e^2} & 0 \\
        -3 e \left(1 + \frac{1}{k_0}\right) \frac{\sin{f_0}}{1-e^2} & 1 & 0 & \frac{e k_0 \cos{f_0} - 2}{1-e^2} & -e \frac{1 + k_0}{1-e^2} \sin{f_0} & 0 \\
        0 & 0 & \sin{f_0} & 0 & 0 & \cos{f_0} \\
        0 & 0 & \cos{f_0} & 0 & 0 & - \sin{f_0}
    \end{bmatrix} \begin{bmatrix}
        \Tilde{\rho} \\
        \theta \\
        \phi \\
        \Tilde{\rho}' \\
        \theta' \\
        \phi'
    \end{bmatrix}_0
\end{equation}
The product of the matrices in Equations~(\ref{eq:YA}) and (\ref{eq:YA_inv}) is the famous YA state transition matrix for relative motion on eccentric orbits.

\section{Second-Order Curvilinear Solution}

Higher-order solutions to the equations of relative motion in spherical coordinates may be found by treating the true solution as a series expansion
\begin{equation}\label{eq:sol_expansion}
    \begin{aligned}
        %\mathbf{\Tilde{\rho}} &= \mathbf{\Tilde{\rho}}_1 + \mathbf{\Tilde{\rho}}_2 + \mathbf{\Tilde{\rho}}_3 + \cdots \\
        \begin{bmatrix}
            \Tilde{\rho} \\ \theta \\ \phi
        \end{bmatrix} &= \begin{bmatrix}
            \Tilde{\rho}_1 \\ \theta_1 \\ \phi_1
        \end{bmatrix} + \begin{bmatrix}
            \Tilde{\rho}_2 \\ \theta_2 \\ \phi_2
        \end{bmatrix} + \begin{bmatrix}
            \Tilde{\rho}_3 \\ \theta_3 \\ \phi_3
        \end{bmatrix} + \cdots
    \end{aligned}
\end{equation}
in which $(\Tilde{\rho}_1,\theta_1,\phi_1)$ capture effects up to $\mathcal{O}(\delta r/r)$, $(\Tilde{\rho}_2,\theta_2,\phi_2)$ capture effects up to $\mathcal{O}(\delta r^2/r^2)$, and so forth. The first-order solution $(\Tilde{\rho}_1,\theta_1,\phi_1)$ to the spherical coordinate dynamics in Equation~(\ref{eq:eom2_curv}) is precisely the YA solution in Equation~(\ref{eq:YA}). For convenience, the initial conditions of $(\Tilde{\rho}_i,\theta_i,\phi_i)$ and $(\Tilde{\rho}'_i,\theta'_i,\phi'_i)$ are chosen to be zero for $i > 1$. The first-order solution is therefore exact at the initial state, i.e. $(\Tilde{\rho}(f_0),\theta(f_0),\phi(f_0))=(\Tilde{\rho}_1(f_0),\theta_1(f_0),\phi_1(f_0))$, and the higher-order components account for the accumulation of error in the first-order solution. This assumption is beneficial because it allows us to use Equation~(\ref{eq:YA_inv}) to define the integration constants $K_1$ through $K_6$ without having to invert a higher-order system.

To derive the second-order solution $(\Tilde{\rho}_2,\theta_2,\phi_2)$, we substitute Equation~(\ref{eq:sol_expansion}) into the equations of motion and expand in products of the components. The only terms in the expansion that contribute to the second-order solution are those that are linear in $(\Tilde{\rho}_2,\theta_2,\phi_2)$ or quadratic in $(\Tilde{\rho}_1,\theta_1,\phi_1)$. Terms involving products of $(\Tilde{\rho}_1,\theta_1,\phi_1)$ and $(\Tilde{\rho}_2,\theta_2,\phi_2)$ components will contribute to the third-order solution and terms quadratic in $(\Tilde{\rho}_2,\theta_2,\phi_2)$ will contribute to the fourth-order solution. Higher-order effects due to terms truncated in the derivation of Equation~(\ref{eq:eom2_curv}) from Equation~(\ref{eq:eom_curv_exact}) will be at least $\mathcal{O}(\delta r^3/r^3)$ and have no contribution to $(\Tilde{\rho}_2,\theta_2,\phi_2)$. Thus, the second-order components solve the system formed by substituting the first-order solution into the nonlinear terms on the right-hand side of Equation~(\ref{eq:eom2_curv}),
\begin{equation}\label{eq:eom2_pert}
    \begin{aligned}
        \Tilde{\rho}_2'' - 2\theta_2' - \frac{3}{k} \Tilde{\rho}_2 &= RHS_\rho (f) =  -\frac{3}{k} \Tilde{\rho}_1^2 + 2 \Tilde{\rho}_1 \theta_1' + \phi_1'^2 + \theta_1'^2 - \phi_1^2 \\
        \theta_2'' + 2\Tilde{\rho}_2' &= RHS_\theta (f) = -2 \Tilde{\rho}_1' \theta_1' + 2 \phi_1' \phi_1 + 2 \Tilde{\rho}_1 \Tilde{\rho}_1' \\
        \phi_2'' + \phi_2 &= RHS_\phi (f) = -2 \theta_1' \phi_1 - 2 \Tilde{\rho}_1' \phi_1'
    \end{aligned}
\end{equation}
where the functions $RHS_i(f)$ have been introduced for generality and to simplify the equations below.

The system in Equation~(\ref{eq:eom2_pert}) simplifies the dynamics of Equation~(\ref{eq:eom2_curv}) by decoupling the out-of-plane component $\phi_2$ from the in-plane components $\Tilde{\rho}_2$ and $\theta_2$. The in-plane equations can be decoupled by integrating the $\theta_2''$ equation once to obtain the system,
\begin{equation}\label{eq:eom2_xy}
    \begin{aligned}
        \Tilde{\rho}_2'' - 2\theta_2' - \frac{3}{k} \Tilde{\rho}_2 &= RHS_\rho (f) \\
        \theta_2' &= -2\Tilde{\rho}_2 + \int RHS_\theta (f) df + c_{\theta1}
    \end{aligned}
\end{equation}
where the constant of integration $c_{\theta 1}$ has been explicitly removed from the integral on the right-hand side. Applying the zero initial conditions to $\theta_2'$ and $\Tilde{\rho}_2$, we find that $c_{\theta1} = -\int RHS_\theta (f) df |_{f_0}$. The integrals in $c_{theta1}$ and Equation~(\ref{eq:eom2_xy}) can be evaluated in terms of the state variables without having to substitute the first-order solution. As an interesting side note, this is different from the rectilinear case. Despite having more terms on the right-hand side of the dynamics in Equation~(\ref{eq:eom2_curv}) than Equation~(\ref{eq:eom2_cart}), the solution derivation is in some ways cleaner in spherical coordinates than in Cartesian. Using integration by parts, we find
\begin{equation}
    \begin{aligned}\label{eq:eom2_y}
    \theta_2' &= -2\Tilde{\rho}_2 + \int RHS_\theta (f) df + c_{\theta1} \\
     &= -2\Tilde{\rho}_2 + \int \left(-2 \Tilde{\rho}_1' \theta_1' + 2 \phi_1' \phi_1 + 2 \Tilde{\rho}_1 \Tilde{\rho}_1'\right) df \\
     &= -2\Tilde{\rho}_2 - 2\theta_1' \Tilde{\rho}_1 + \phi_1^2 - \Tilde{\rho}_1^2 + c_{\theta 1}
    \end{aligned}
\end{equation}
where $c_\theta 1 = 2\theta'(f_0) \Tilde{\rho}(f_0) - \phi(f_0)^2 + \Tilde{\rho}(f_0)^2$.
Substituting the expression for $\theta_2'$ into the equation for $\Tilde{\rho}_2''$ leads to the second-order linear inhomogeneous ODE,
\begin{equation}\label{eq:eom2_x}
    \Tilde{\rho}_2'' + \left(4 - \frac{3}{k}\right) \Tilde{\rho}_2 = -\left(2 + \frac{3}{k}\right)\Tilde{\rho}_1^2 - 2\Tilde{\rho}_1 \theta_1' + \phi_1'^2 + \theta_1'^2 + \phi_1^2 + 2 c_{\theta1}
\end{equation}
Second-order components of the relative motion appear only on the left of Equation~(\ref{eq:eom2_x}), while the right may be written as an explicit function of $f$ using Equation~(\ref{eq:YA}).

Equation~(\ref{eq:eom2_x}) can be solved by variation of parameters if two linearly independent solutions are available for the homogeneous equation\cite{BoyceDiPrima}
\begin{equation}\label{eq:x_homog}
    \Tilde{\rho}''_2 + \left( 4 - \frac{3}{k} \right) \Tilde{\rho}_2 = 0
\end{equation}
Because the higher-order terms involving $\Tilde{\rho}_1$ and $\theta_1$ do not appear in the homogeneous equation, it is identical to that obtained from the TH equations. The solutions to this equation introduced by Yamanaka and Ankersen are
\begin{equation}\label{eq:YA_xsols}
    \begin{aligned}
        \varphi_1 &= k \sin{f} \\
        \varphi_2 &= 3 e^2 k J(t) \sin{f} + k \cos{f} - 2 e
    \end{aligned}
\end{equation}
and their linear independence was demonstrated in that work.\cite{YA} The particular solution $\varphi_p$ to any inhomogeneous equation formed by placing an arbitrary function $RHS(f)$ of the independent variable on the right of Equation~(\ref{eq:x_homog}) can be found using the variation of parameters formula,
\begin{equation}\label{eq:vop}
    \varphi_p = \varphi_1 \int \frac{\varphi_2 RHS(f)}{1-e^2} df - \varphi_2 \int \frac{\varphi_1 RHS(f)}{1-e^2} df
\end{equation}
where the Wronskian in the denominator is $\varphi_1' \varphi_2 - \varphi_1 \varphi_2' = 1-e^2$. By superposition, the general solution is the sum of the particular solution and a linear combination of the homogeneous solutions $\varphi_1$ and $\varphi_2$:
\begin{equation}\label{eq:xh+xp}
    \Tilde{\rho}_2 = c_{\rho1} \varphi_1 + c_{\rho2} \varphi_2 + \varphi_p
\end{equation}
The particular solution $\varphi_p$ is found by combining Equations~(\ref{eq:YA}), (\ref{eq:eom2_x}), (\ref{eq:YA_xsols}), and (\ref{eq:vop}). The constants $c_{\rho1}$ and $c_{\rho2}$ are found by satisfying the zero initial conditions, $\Tilde{\rho}_2(f_0) = 0$ and $\Tilde{\rho}_2'(f_0) = 0$.

Having solved for $\Tilde{\rho}_2$, the along-track correction $\theta_2$ may be found by direct integration of Equation~(\ref{eq:eom2_y}), along with the zero initial conditions. Finally, applying the variation of parameters procedure to the third line of Equation~(\ref{eq:eom2_curv}) using the homogeneous solutions $\varphi_1 = \sin{f}$ and $\varphi_2 = \cos{f}$ and $\varphi_1' \varphi_2 - \varphi_1 \varphi_2' = 1$ results in the out-of-plane correction $\phi_2$. 

Combining the expressions for $\Tilde{\rho}_2$, $\theta_2$, and $\phi_2$ with the first-order components $\Tilde{\rho}_1$, $\theta_1$, and $\phi_1$ from Equation~(\ref{eq:YA}), we obtain the new solution to the curvilinear equations of relative motion accurate to second-order in the normalized coordinates,
\begin{equation}\label{eq:eom2_sol}
    \begin{aligned}
    \Tilde{\rho} &= K_1 \left( 1 - \frac{3}{2} e k J(t) \sin{f} \right) + K_2 k \sin{f} + K_3 k \cos{f} \\
     & \quad + c_{\rho j} \left(1 - \frac{3}{2} e k J(t) \sin{f} \right) + c_{\rho s} k \sin{f} + c_{\rho c} k \cos{f} \\
     & \quad + K_1^2 \left(\frac{1}{4} + \frac{9}{8} k^3 J(t)^2 e \cos f \right) - \frac{3}{2} \left( K_1 K_2 \cos f - K_1 K_3 \sin f \right) k^3 J(t) \\
     & \quad + K_2^2 \left[ \left(-\frac{e^2}{2} \sin^2 f + \frac{3}{2}(k-1) + \frac{1}{1-e^2} \right)\cos^2 f + \frac{e(1+e^2) \cos f}{2(1-e^2)}\right] \\
     & \quad + K_2 K_3 \frac{\left( e k^2 - (1+k)\cos f\right)k \sin f}{1-e^2} + K_3^2 \frac{k(3-k-k^2+k^3 - (1+k)(e^2+\cos^2 f)}{2(1-e^2)} \\
     \\[-2pt]
     \theta &= K_4 + K_2 (1 + k) \cos{f} - K_3 (1 + k) \sin{f} - \frac{3}{2} K_1 k^2 J(t) \\
     & \quad + \left(c_{\rho s} - K_1 K_2\right) \left( (1 + k) \cos{f} - (1+k_0) \cos{f_0} \right) - \frac{3}{2} \left(K_1^2 - K_1 K_3 e - c_{\rho j} \right) k^2 J(t) \\
     & \quad + \left(K_1 K_3 - K_2^2 \frac{e^3}{2(1-e^2)} - c_{\rho c} \right)\left((1+k) \sin{f} - (1+k_0)\sin{f_0} \right) \\
     & \quad + K_1^2 \left(-\frac{9}{4} e k^3 J(t)^2 \sin f \right) + 3 \left(K_1 K_2 \sin f + K_1 K_3 \cos f \right) k^3 J(t) \\
     & \quad + \left(K_3^2 - K_2^2 \right) \left[ \left(\frac{\cos f + 2e}{2(1-e^2)} + k(1+k) \cos f \right) \sin f \right. \\
     & \qquad \qquad \qquad \qquad \left. - \left(\frac{\cos f_0 + 2e}{2(1-e^2)} + k_0(1+k_0) \cos f_0 \right)\sin f_0 \right] \\
     & \quad + K_2 K_3 \left[ \left( k^2 + \frac{k^2}{1-e^2} - (1 + 2k + 2k^2) \cos^2 f \right) \right. \\
     & \qquad \qquad \qquad \left. - \left(k_0^2 + \frac{k_0^2}{1-e^2} - (1 + 2k_0 + 2k_0^2) \cos^2 f_0 \right) \right] \\
     & \quad + K_3^2 e \left( \sin f - \sin f_0 \right) + \frac{1}{4}\left( K_6^2 - K_5^2 \right)\left(\sin 2f - \sin 2f_0 \right) + K_5 K_6 \left( \sin^2 f - \sin^2 f_0 \right) \\[0pt]
     \\[-2pt]
    \phi &= K_5 \sin{f} + K_6 \cos{f} \\
     & \quad + \frac{3}{2}\left( K_1 K_6 \sin f - K_1 K_5 \cos f \right)k^2 J(t) + \frac{3}{2} \left(  K_1 K_5 \cos f_0 -  K_1 K_6 \sin f_0 \right) \sin(f-f_0) \\
     & \quad + 2\left((K_2 K_5 - K_3 K_6) \cos f_0 - (K_2 K_6 + K_3 K_5) \sin f_0 \right) k_0 \sin f_0 \sin(f-f_0) \\
     & \quad + K_2 K_5 \left( (1+k) \cos f - (1+k_0) \cos f_0 \right) \cos f \\
     & \quad - \left(K_2 K_6 + K_3 K_5\right) \left( (1+k) \cos f - (1+k_0) \cos f_0 \right) \sin f \\
     & \quad + K_3 K_6 \left( (1+k) \sin^2 f - e \sin^2 f_0 \cos f - 2\sin f_0 \sin f \right)
    \end{aligned}
\end{equation}
where the constants $c_{\rho j}$, $c_{\rho s}$, and $c_{\rho c}$ are provided in the appendix as functions of the integration constants $K_1$ through $K_6$ and the initial true anomaly. As with the YA solution, the integration constants are found from the initial conditions using Equation~(\ref{eq:YA_inv}). 

Although the right-hand side of Equation~(\ref{eq:eom2_pert}) appears more complicated than the second-order Cartesian system in Equation~(\ref{eq:eom2_pert}), these equations and their solution are simpler in terms of the combinations of the integration constants with nonzero coefficients. There are 21 possible pairings of the six constants $\{K_i\}$. The Cartesian solution employs 19 of these--all except $K_4 K_5$ and $K_4 K_6$. The curvilinear solution involves only 15 pairs, having no terms with $K_4$. This constant represents an initial offset in $\theta$ and is absent because only the derivative of $\theta_1$ appears in the equations of motion. Significantly, its absence makes the equations invariant under $\theta$ rotations and allows the solution to retain its accuracy in the presence of large along-track separations.

Finally, consider the limit as $e \rightarrow 0$. This allows us to make the substitutions $k \rightarrow 1$, $J(t) \rightarrow n(t-t_0)$, $f \rightarrow n(t-t_{ref})$, and $f_0 \rightarrow n(t_0 - t_{ref})$. With time normalized such that $n = 1$ and letting $t_0 = t_{ref} = 0$, Equation~(\ref{eq:eom2_sol}) becomes
\begin{equation}
    \begin{aligned}
    \rho &= K_1 + K_2 \sin{t} + K_3 \cos{t} \\
    & \quad - \frac{3}{2} K_1 K_2 t \cos t + \frac{3}{2}K_1 K_3 t \sin t + \frac{1}{2} \left( K_2^2 - K_3^2 \right)(\cos 2 t - 1) - K_2 K_3 \sin 2t \\
    & \quad + \left(\frac{15}{4} K_1^2 + 10 K_1 K_3 - 2K_2^2 + 5 K_3^2 - K_5^2 + K_6^2 \right) (\cos t - 1) + \left(\frac{3}{2} K_1 K_2 + 2 K_2 K_3 \right) \sin t \\
    \\
    \theta &= K_4 + 2 K_2 \cos{t} - 2 K_3 \sin{t} - \frac{3}{2} K_1 t \\
    & \quad  + \left(\frac{15}{2} (K_1^2 + 2 K_1 K_3 + K_3^2) - \frac{3}{2} (K_2^2 + K_5^2 - K_6^2)  \right) t \\
    & \quad + 3 K_1 K_2 t \sin t + 3 K_1 K_3 t \cos t + (K_1 K_2 + 4 K_2 K_3) (\cos t - 1) \\
    & \quad + \left( -\frac{15}{2}K_1^2 - 18 K_1 K_3 + 4 K_2^2 - 10 K_3^2 + 2K_5^2 - 2 K_6^2 \right) \sin t \\
    & \quad + \frac{1}{4}(5K_3^2 - 5K_2^2 + K_6^2 - K_5^2) \sin 2t - \frac{1}{2}(5 K_2 K_3 + K_5 K_6) (\cos 2t - 1) \\ 
    \\
    \phi &= K_5 \sin{t} + K_6 \cos{t} \\
    & \quad + K_2 K_5 + K_3 K_6 + \frac{3}{2} K_1 K_6 t \sin t - \frac{3}{2}K_1 K_5 t \cos t + \left( \frac{3}{2} K_1 K_5 + 2 K_2 K_6 + 2 K_3 K_5 \right) \sin t \\
    & \quad - 2K_2 K_5 \cos t - (K_2 K_6 + K_3 K_5) \sin 2t + (K_2 K_5 - K_3 K_6) \cos 2t \\
    \end{aligned}
\end{equation}
which is identical to the spherical QV solution obtained by extension of CW to second-order.

\section{Validation}

In this section we will compare the performance of the curvilinear solution with its Cartesian cousin as well as other translational state solutions from the literature. These models are evaluated against an unperturbed Keplerian truth to show how well each captures the relative motion subject to the assumptions under which it was derived. This choice reduces the number of parameters needed to fully specify the chief's absolute motion to three: semimajor axis, eccentricity, and true anomaly. Combined with the six parameters needed to characterize the relative motion, we have a nine-dimensional state space. Rather than attempting a full parameter sweep, we will focus on the effects of eccentricity and inter-spacecraft separation.

To maintain consistency and feasibility, all scenarios are initialized at perigee with an altitude $h_p$ of 750~km and propagated for 10 orbits. Because the semimajor axis is larger for more eccentric orbits, the duration of the simulated scenarios can differ. Table~\ref{tab:ICs_abs} summarizes the chief's absolute orbit parameters common to all simulations. 

\begin{table}[t]
    \centering
    \caption{Chief Orbit Parameters for Performance Comparison Scenarios}
    \begin{tabular}{|ccccc|}
        \hline
        $h_p$ (km) & $i$ & $\Omega$ & $\omega$ & $f_0$ \\
        \hline
        750 & 98$^{\circ}$ & 30$^{\circ}$ & 30$^{\circ}$ & 0$^{\circ}$ \\
        \hline
    \end{tabular}
    \label{tab:ICs_abs}
\end{table}

Many alternative representations would serve to specify the relative motion. We will use the quasi-nonsingular relative orbital elements (ROE), which are defined in terms of the Keplerian orbital elements of the chief and deputy as
\begin{equation}\label{eq:ROE}
    \boldsymbol{\updelta \upalpha} = \begin{bmatrix}
        \delta a \\
        \delta \lambda \\
        \delta e_x \\
        \delta e_y \\
        \delta i_x \\
        \delta i_y
    \end{bmatrix} = \begin{bmatrix}
        \frac{a_d - a}{a} \\
         (u_d - u) + (\Omega_d - \Omega) \cos{i} \\
         e_d \cos{\omega_d} - e \cos{\omega} \\
         e_d \sin{\omega_d} - e \sin{\omega} \\
         i_d - i \\
         (\Omega_d - \Omega) \sin{i}
    \end{bmatrix}
\end{equation}
where $u = f + \omega$ is the argument of latitude and subscripts are omitted on elements related to the chief's orbit. Unlike the components of relative position and velocity, all of the ROE are constant for unperturbed orbital motion except $\delta\lambda$. %They also offer direct insight into the relative orbit geometry.  
For near-circular orbits, the extent of the relative motion in the radial direction is proportional to the $L_2$-norm of the relative eccentricity vector $\boldsymbol{\updelta}\mathbf{e} = [\delta e_x, \delta e_y]^T$ and the extent of the out-of-plane motion is proportional to the $L_2$-norm of the relative inclination vector $\boldsymbol{\updelta}\mathbf{i} = [\delta i_x, \delta i_y]^T$. The mean along-track separation is given by $a \delta \lambda$ and along-track drift is governed by $a \delta a$.

\begin{figure}[htb]
	\centering\includegraphics[width=5.5in,trim=0.2in 0in 0.5in 0in, clip]{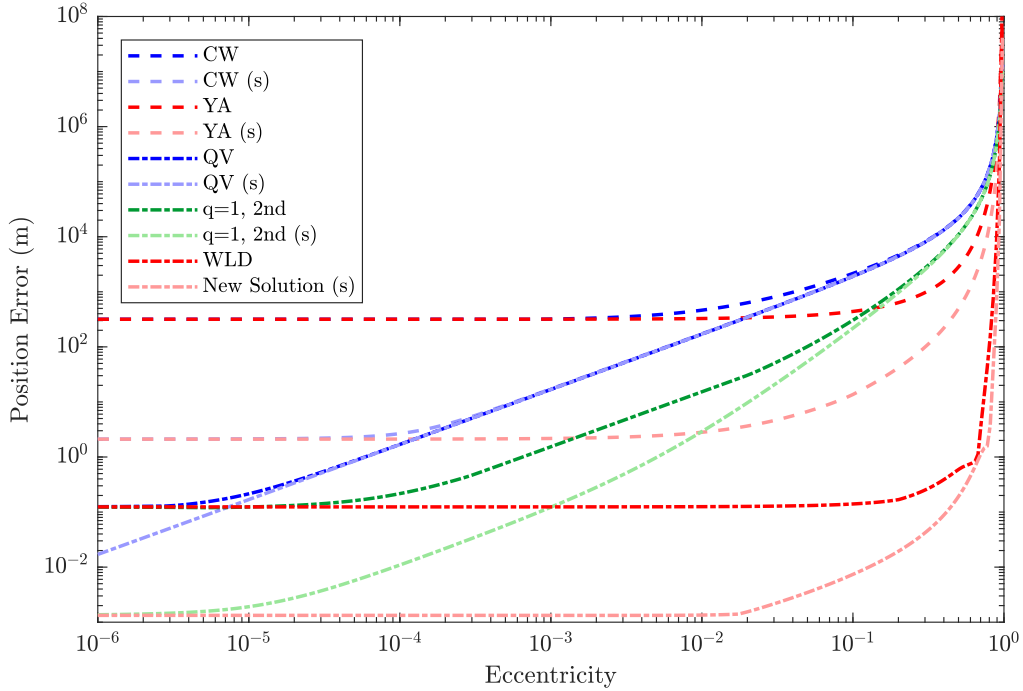}
	\caption{Maximum Position Error Against Eccentricity with \textit{\textbf{a\!}} $\mathbf{\boldsymbol{\delta \alpha} =}$ [0, 0, 0, 2, 0, 2] km.}
	\label{fig:ecc_000101}
\end{figure}

Figure~\ref{fig:ecc_000101} compares the maximum error over 10 orbits of several relative motion solutions against the eccentricity of the chief's orbit, with the ROE $a \boldsymbol{\updelta \upalpha} = [0, 0, 0, 2, 0, 2]^T$~km. The comparison includes both Cartesian and curvilinear versions of CW, YA, QV, and a solution that treats the chief's eccentricity as a perturbation to the circular orbit dynamics. The latter is characterized by the parameter $q$ that represents the perturbing strength of eccentricity relative to inter-spacecraft separation, as well as the maximum order in separation of eccentricity terms included in the solution. For $q=1$, the solution includes terms up to $e \delta r$ and $\delta r^2$ and is therefore second order in separation.\footnote{A small correction to the spherical coordinate equations of motion for slightly eccentric orbits was necessary, and is described in the Appendix} All solutions are style- and color-coded according to their order and underlying assumptions, respectively. Dashed lines indicate a linear model and dash-dot a second-order model. Blue lines are used for models that assume circular orbits, red fully incorporate eccentricity through coordinate transformation, and green treat eccentricity as a perturbation to the circular orbit dynamics. Curvilinear solutions are plotted using lighter shades and marked with (s) in the legend to emphasize the use of spherical coordinates. This formatting matches that used in previous comparisons by the authors, though several higher-order solutions included in earlier work are omitted here for clarity.\cite{WillisLovellDAmico}

The log-log representation of Figure~\ref{fig:ecc_000101} reveals the grouping of solutions according to their order, coordinate system, and underlying assumptions. For very low eccentricities, the linear CW and YA solutions give the highest error, with the curvilinear models giving a 100-fold improvement in accuracy over the rectilinear models. The rectilinear second-order solutions give another order of magnitude improvement, and the curvilinear three orders of magnitude. Those solutions that assume circular orbits or treat eccentricity as a perturbation diverge from the YA and new solution as eccentricity increases. In this family of scenarios, the new solution in spherical coordinates is the most accurate by several orders of magnitude for eccentricities between 0.0001 and 0.9. The presence of $1-e^2$ in the denominator of terms in Equations~(\ref{eq:YA_inv}) and (\ref{eq:eom2_sol}) makes these solutions singular for parabolic orbits and causes the error to diverge near $e = 1$.

\begin{figure}[htb]
	\centering\includegraphics[width=5.5in,trim=0.2in 0in 0.5in 0in, clip]{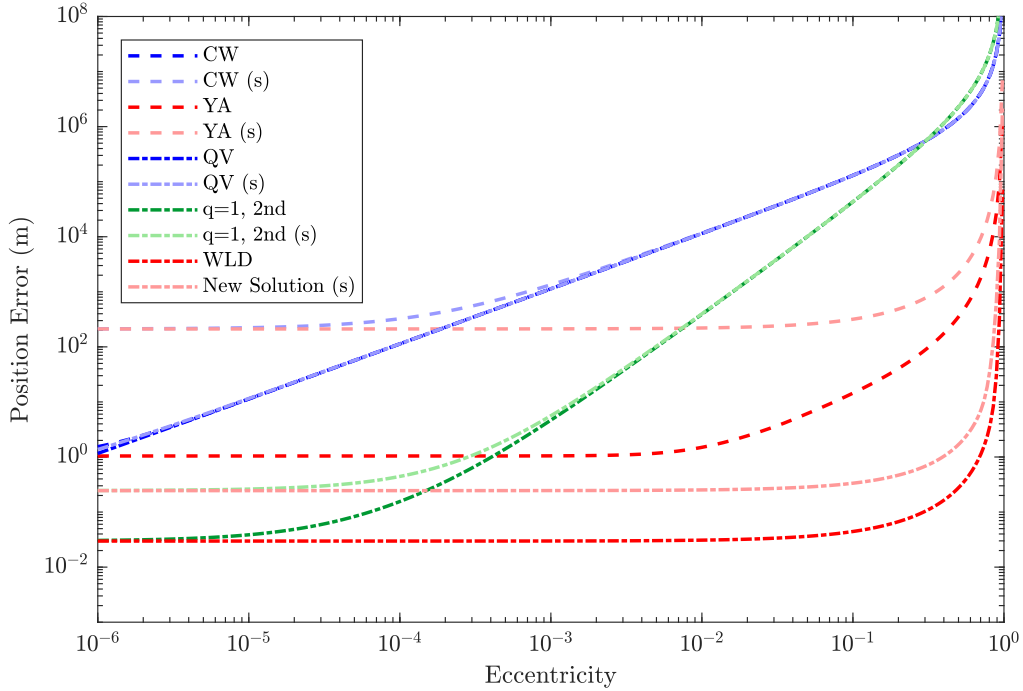}
	\caption{Maximum Position Error Against Eccentricity with \textit{\textbf{a\!}} $\mathbf{\boldsymbol{\delta \alpha} =}$ [0, 0, 2, 0, 2, 0] km.}
	\label{fig:ecc_001010}
\end{figure}

When comparing rectilinear and curvilinear models, it is important to examine the sensitivities of each to the initial conditions. Figure~\ref{fig:ecc_001010} gives the same comparison as Figure~\ref{fig:ecc_000101}, but with the ROE $a \boldsymbol{\updelta \upalpha} = [0, 0, 2, 0, 2, 0]^T$~km. This rotation of the relative eccentricity and inclination vectors has no change on the shape of the relative motion, but alters its initial phase by 90\degree. In the previous scenario, the deputy was initialized at its maximum along-track and cross-track separation. Now, the initial separation is purely radial. This has a dramatic effect on the propagation errors and especially on which solutions are more accurate. When the initial separation is radial, the rectilinear models outperform the curvilinear models. The error is dominated by along-track drift, brought on by imperfect modeling of the relative semimajor axis. The rectilinear and curvilinear models differ in the phase at which the models most accurately describe the relative motion. Incorporating higher-order terms reduces these modeling errors and reduces sensitivity to the initial phase. As a result, the second-order solutions in both coordinate systems are more accurate than the linear solutions in either scenario.

\begin{figure}[htb]
	\centering\includegraphics[width=5.3in,trim=0.2in 0in 0.55in 0in, clip]{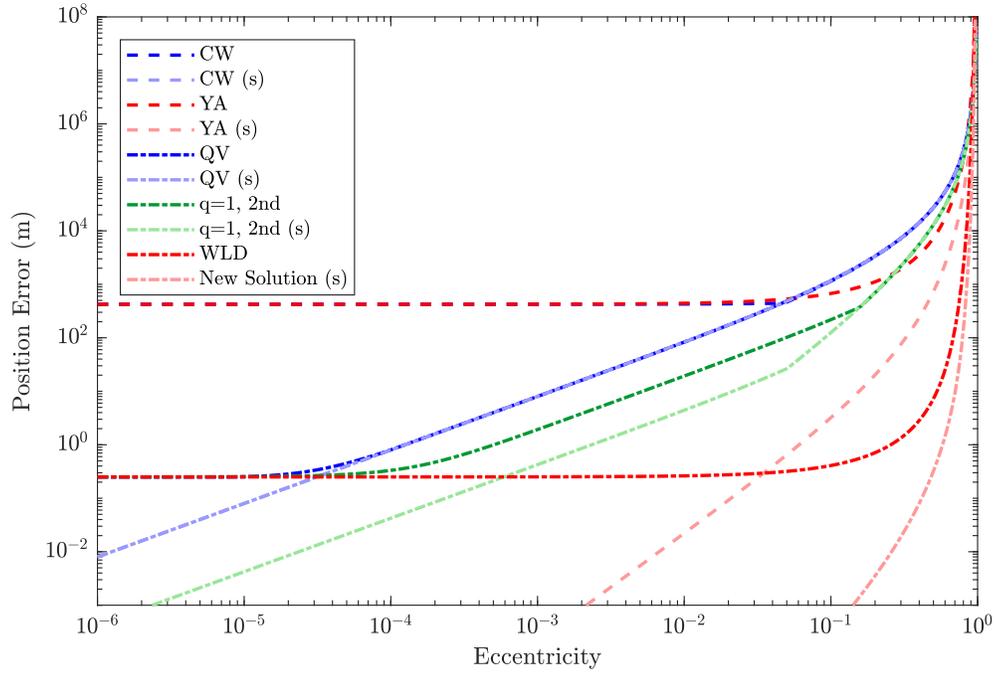}
	\caption{Maximum Position Error Against Eccentricity with \textbf{\textit{a\!}} $\mathbf{\boldsymbol{\delta \alpha} =}$ [0, 4, 0, 0, 0, 0] km.}
	\label{fig:ecc_010000}
\end{figure}

The scenarios in Figures~\ref{fig:ecc_000101} and \ref{fig:ecc_001010} involve centered relative motion. Figure~\ref{fig:ecc_010000} compares the solutions for a case where both spacecraft lie on the same orbit, with a separation in true anomaly represented by the single nonzero ROE $a \delta \lambda = 4$~km. For near-circular orbits, the curvilinear models approach arbitrary levels of accuracy. As noted in the derivation of the new solution, the second-order equations of motion do not depend on the along-track separation, allowing for very high accuracy propagation in the presence of a constant phase offset.

\begin{figure}[t!]
	\centering\includegraphics[width=5.3in,trim=0.2in 0in 0.55in 0in, clip]{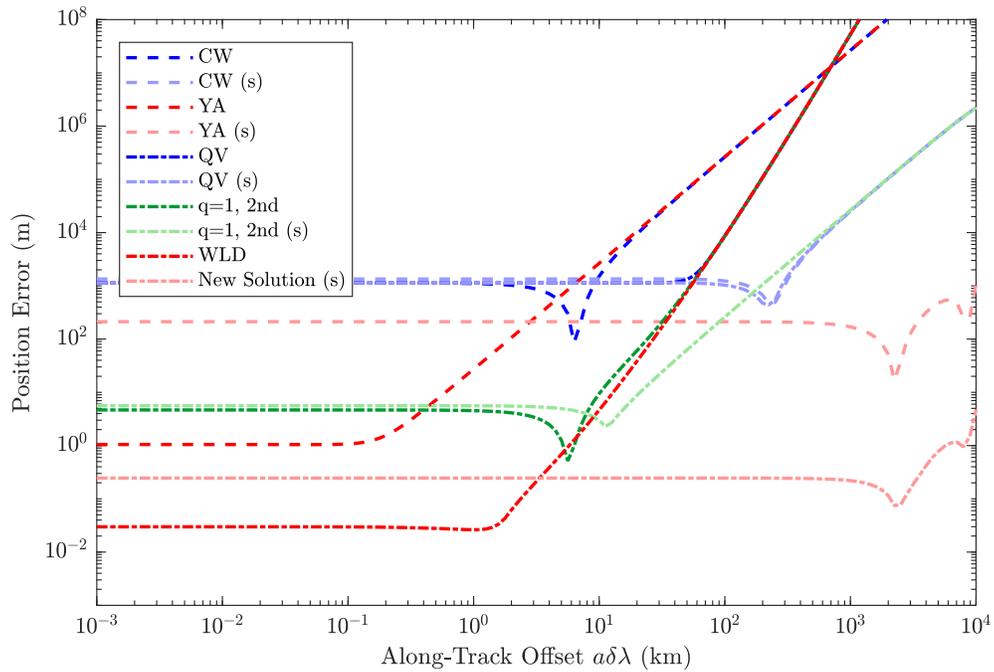}
	\caption{Maximum Position Error Against Along-Track Separation.}
	\label{fig:sep_011010_001}
\end{figure}

The performance of the solutions against increasing separation is as important as that against eccentricity. Figure~\ref{fig:sep_011010_001} compares the maximum position error of the same set of translational state solutions against a Keplerian truth model over 10 orbits as a function of separation. The scenario uses $e = 0.001$, fixed relative eccentricity and inclination vectors, and a range of along-track offsets. The ROE are $a \boldsymbol{\updelta \upalpha} = [0, a \delta \lambda, 2, 0, 2, 0]^T$~km. The phases of the relative eccentricity and inclination vectors have been chosen to match the case of Figure~\ref{fig:ecc_001010}, in which the new second-order solution was less accurate than its rectilinear counterpart. Indeed, a vertical section from the left side of Figure~\ref{fig:sep_011010_001} would match a slice of Figure~\ref{fig:ecc_001010} along $e = 0.001$. However, that is for centered relative motion. The Cartesian solution loses accuracy for along-track offsets greater than 1~km, which would still place the chief within the in-plane projection of the relative motion. The new, second-order curvilinear solution does not lose accuracy until the along-track offset is more than 1000 km, and has meter-level accuracy well beyond the horizon in this scenario at 750~km altitude. 

A striking feature of Figure~\ref{fig:sep_011010_001} is the sharp drop in position error experienced by most solutions, both in Cartesian and spherical coordinates. Recall that the position error is dominated by the along-track drift which results from imperfect representation of the energy-matching (no-drift) condition in the approximate solutions. These sharp accuracy improvements occur where the the direction of the along-track drift reverses. With our initialization strategy, changing the along-track offset has a small effect on the initial phase of the relative motion and thus how well the no-drift condition is captured by the different models. The propagation errors in the spherical coordinate solutions for eccentric orbits are small enough for two such accuracy spikes to appear in the range of offsets shown. The other models experience only one spike, if any, before other effects dominate the propagation error.

Up to this point we have exclusively considered translational state solutions that are closely related to the new solution. However, many authors favor solutions based on orbital elements, either through orbit element differences or the ROE defined in Equation~(\ref{eq:ROE}). Because the orbital elements are constants of motion in the two-body problem, they map to relative position and velocity with zero error compared to a Keplerian truth. For control system design in the orbital element state space, it is preferable to approximate the dynamics so that the solution is linear in the state variables. In the case of ROE, this affects only the mean relative argument of latitude $\delta \lambda$ because each of the other ROE are constant with respect to Keplerian dynamics. Expanding the nonlinear terms to second order, $\delta \lambda$ is given by
\begin{equation}\label{eq:roe_approx}
\delta \lambda \approx \delta \lambda_0 - \frac{3}{2} \delta a \, nt + \frac{15}{8} \delta a^2 \frac{n t}{a} 
\end{equation}
where $n$ is the mean motion of the chief spacecraft.

\begin{figure}[b!]
	\centering\includegraphics[width=5.5in,trim=0.1in 0in 0.2in 0in, clip]{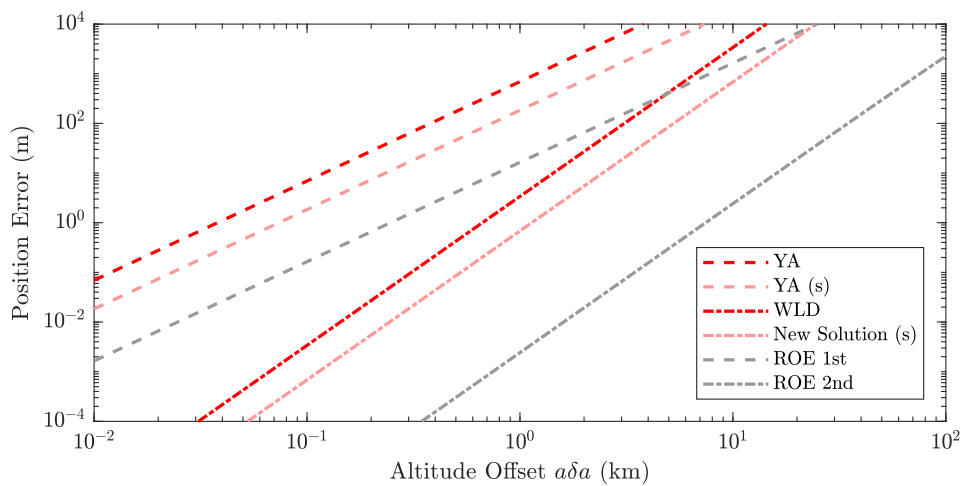}
	\caption{Maximum Position Error Against Altitude Offset.}
	\label{fig:roe_test}
\end{figure}

Figure~\ref{fig:roe_test} compares the performance of YA and the second-order eccentric solutions in rectilinear and curvilinear coordinates with that obtained by propagating the ROE using the approximate model in Equation~(\ref{eq:roe_approx}). The scenario chosen for this test has a difference in semimajor axis only and a moderately large eccentricity of 0.1. The figure shows the maximum position error over ten orbits and includes solutions for $\delta \lambda$ truncated at first and second order in $\delta a$ and uses the exact nonlinear mapping from ROE to relative position and velocity. This test is of particular interest because it involves an initially radial offset that grows into a large along-track error. The first-order ROE model is more accurate than the spherical YA solution, which is more accurate than the cartesian YA solution. However, the second-order translational state solutions are more accurate when $a \delta a$ is no more than a few kilometers in this LEO scenario. The second-order ROE propagation is far more accurate than the translational state solutions, and the first-order propagation is more accurate for sufficiently large $a \delta a$ because the error grows more rapidly for the higher-order solutions. Note that the order of both the ROE and translational state models is reflected in the slope of the error trends on this log-log plot.

\section{Conclusion}

A new, second-order solution for the relative motion of two spacecraft on eccentric orbits has been derived in spherical coordinates. It is related to the second-order Cartesian solution previously derived by the authors through its connection to the Yamanaka-Ankersen state transition matrix. However, the two second-order solutions are not equivalent because the new solution is not directly based on relative position and velocity and the higher-order effects differ in rectilinear and curvilinear coordinates. The new solution was validated against an unperturbed Keplerian truth model and compared with several related solutions from the literature. In all test cases, the new solution gave a thousand-fold improvement in accuracy over the first-order curvilinear solution. Whether the second-order solution in spherical coordinates or its Cartesian counterpart is more accurate depends in part upon the initial conditions, but the spherical solution is generally better in the presence of large along-track separations. 

Both the second-order rectilinear and curvilinear solutions assume unperturbed, Keplerian motion. However, perturbations from Earth oblateness, solar radiation pressure, atmospheric drag, and third body tidal effects can have as large an effect on the relative motion as the higher-order Keplerian dynamics. Future work on this topic should explore the solution's sensitivity to uncertainty in the absolute and relative states and accuracy in the presence of perturbing forces. Going a step further, the same methodology employed to develop the second-order corrections may be used to incorporate the leading-order effects of such disturbances into the solution.

%Both the second-order rectilinear and curvilinear solutions assume unperturbed, Keplerian motion. A framework has been presented for extending these solutions to include general perturbations and the equations governing the corrections for $J_2$ perturbation were derived. Perturbation effects manifest in the relative dynamics in four distinct ways: through the differential inertial force acting on the formation, through the perturbation of the rotating reference frame in which the motion is described, through the transformation of the independent variable from time to argument of latitude, and finally through the effect of the osculating orbital elements on the parameters appearing in the leading-order dynamics. It has been shown through numerical integration that an order of magnitude improvement in accuracy can be achieved if all effects are taken into account.

%The perturbation framework was derived in the context of the rectilinear dynamics model. However, because it is linear in the relative position and velocity coordinates, it is believed that it will have the same form in curvilinear coordinates. Future work on this topic should demonstrate this connection, solve the set of equations established for $J_2$ perturbation, and develop solutions for other perturbations of interest. 

% Assess sensitivity of equations to perturbations.

%\section{Acknowledgment}
%Any acknowledgments by the author may appear here. The acknowledgments section is optional.

%\section{Notation}
%\begin{tabular}{r l}
%	$a$ & a real number \\
%	$b$ &  the square root of $a$ \\
%\end{tabular} \\

\section{ACKNOWLEDGMENTS}

This work was supported by a NASA Space Technology Research Fellowship. The authors would also like to thank Eric Butcher and Ethan Burnett for providing the code for their higher-order solutions.

\bibliographystyle{AAS_publication}   % Number the references.
\bibliography{references}   % Use references.bib to resolve the labels.

\begin{thebibliography}{10}

\bibitem{DistributedSS}
S.~D'Amico, M.~Pavone, S.~Saraf, A.~Alhussien, T.~Al-Saud, S.~Buchman,
  R.~Bryer, and C.~Farhat, ``Distributed Space Systems for Future Science and
  Exploration,''  8th International Workshop on Spacecraft Formation Flying,
  Delft University, June 8-10, 2015.

\bibitem{SullivanGrimberg}
J.~Sullivan, S.~Grimberg, and S.~D'Amico, ``Comprehensive Survey and Assessment
  of Spacecraft Relative Motion Dynamics Models,''  {\em Journal of Guidance,
  Control, and Dynamics}, Vol.~40, No.~8, 2017, pp.~1837--1859.

\bibitem{HCW}
W.~H. Clohessy and R.~S. Wiltshire, ``Terminal Guidance System for Satellite
  Rendezvous,''  {\em Journal of Guidance, Control, and Dynamics}, Vol.~27,
  No.~9, 1960, pp.~653--658.

\bibitem{London}
H.~S. London, ``Second Approximation to the Solution of the Rendezvous
  Equations,''  {\em AIAA Journal}, Vol.~1, No.~7, 1963, pp.~1691--1693.

\bibitem{SasakiQV}
M.~L. Anthony and F.~T. Sasaki, ``Rendezvous Problem for Nearly Circular
  Orbits,''  {\em AIAA Journal}, Vol.~3, No.~7, 1965, pp.~1666--1673.

\bibitem{NewmanQV}
M.~T. Stringer, B.~A. Newman, T.~A. Lovell, and A.~Omran, ``Analysis of a New
  Nonlinear Solution of Relative Orbital Motion,''  23rd International
  Symposium on Space Flight Dynamics, Pasadena, California, October 29-November
  2 2012.

\bibitem{QVcompare}
B.~A. Newman, A.~J. Sinclair, T.~A. Lovell, and A.~Perez, ``Comparison of
  Nonlinear Analytical Solutions for Relative Orbital Motion,''  AIAA/AAS
  Astrodynamics Specialist Conference, San Diego, California, August 4-7 2014.

\bibitem{Melton}
R.~G. Melton, ``Time-Explicit Representation of Relative Motion Between
  Elliptical Orbits,''  {\em Journal of Guidance, Control, and Dynamics},
  Vol.~23, No.~4, 2000, pp.~604--610.

\bibitem{ButcherLovellCart}
E.~A. Butcher, T.~A. Lovell, and A.~Harris, ``Third Order Cartesian Relative
  Motion Perturbation Solutions for Slightly Eccentric Chief Orbits,''  26th
  AAS/AIAA Space Flight Mechanics Meeting, Napa, CA, February 14-18, 2016.

\bibitem{TH}
J.~Tschauner and P.~Hempel, ``Optimale Beschleunigungsprogramme fur das
  Rendezvous-Manover,''  {\em Astronautica Acta}, Vol.~10, No.~5-6, 1964,
  p.~296.

\bibitem{CarterTH}
T.~E. Carter, ``State Transition Matrices for Terminal Rendezvous Studies:
  Brief Survey and New Example,''  {\em Journal of Guidance, Control, and
  Dynamics}, Vol.~21, No.~1, 1998, pp.~148--155.

\bibitem{YA}
K.~Yamanaka and F.~Ankersen, ``New State Transition Matrix for Relative Motion
  on an Arbitrary Elliptical Orbit,''  {\em Journal of Guidance, Control, and
  Dynamics}, Vol.~25, No.~1, 2002, pp.~60--66.

\bibitem{WillisLovellDAmico}
M.~Willis, A.~Lovell, and S.~D'Amico, ``Second Order Analytical Solution for
  Relative Motion on Arbitrarily Eccentric orbits,''  AIAA/AAS Space Flight
  Mechanics Meeting, Ka'anapali, Maui, HI, January 13-17 2019.

\bibitem{ScheeresSTT}
R.~S. Park and D.~J. Scheeres, ``Nonlinear Mapping of Gaussian Statistics:
  Theory and Applications to Spacecraft Trajectory Design,''  {\em Journal of
  Guidance, Control, and Dynamics}, Vol.~29, No.~6, 2006, pp.~1367--1375.

\bibitem{AlfriendFF}
K.~T. Alfriend, S.~R. Vadali, P.~Gurfil, J.~P. How, and L.~S. Breger, {\em
  Spacecraft Formation Flying: Dynamics, control and navigation}.
\newblock Elsevier, 2010.

\bibitem{ButcherBurnettCurv}
E.~A. Butcher, E.~Burnett, and T.~A. Lovell, ``Comparison of Relative Orbital
  Motion Perturbation Solutions in Cartesian and Spherical Coordinates,''  27th
  AAS/AIAA Space Flight Mechanics Meeting, San Antonio, TX, February 5-9, 2017.

\bibitem{ButcherLovellDT}
E.~A. Butcher and T.~A. Lovell, ``Spherical Coordinate Perturbation Solutions
  to Relative Motion Equations: Application to Double Transformation Spherical
  Solution,''  26th AAS/AIAA Space Flight Mechanics Meeting, Napa, CA, February
  14-18, 2016.

\bibitem{HanYASpherical}
C.~Han, H.~Chen, and G.~A. e.~al., ``A linear model for relative motion in an
  elliptical orbit based on a spherical coordinate system,''  {\em Acta
  Astronautica}, Vol.~157, April 2019, pp.~465--476.

\bibitem{BoyceDiPrima}
W.~E. Boyce and R.~C. DiPrima, {\em Elementary Differential Equations}.
\newblock Wiley, 2008.

\end{thebibliography}

\newpage
\appendix
\section*{APPENDIX: Coordinate Transformations}
\subsection*{Relative Position and Velocity to Spherical Coordinates}
To convert relative position $\mathbf{\boldsymbol \updelta r} = [x,y,z]^T$ and relative velocity $\mathbf{\boldsymbol \updelta v} = [\dot{x},\dot{y},\dot{z}]^T$ vectors to spherical coordinates, use the transformations 
\begin{equation}
    \begin{aligned}
        \rho &= \sqrt{(r + x)^2 + y^2 + z^2} - r \\
        \theta &= \tan^{-1} \left(\frac{y}{r+x} \right) \\
        \phi &= \sin^{-1} \left( \frac{z}{r+\rho} \right) \\
        \dot{\rho} &= \frac{(r+x)(\dot{r}+\dot{x})+y\dot{y}+z\dot{z}}{r+\rho} - \dot{r} \\
        \dot{\theta} &= \frac{(r+x)\dot{y} - y (\dot{r}+\dot{x})}{(r+x)^2 + y^2} \\
        \dot{\phi} &= \frac{(r+\rho)\dot{z} - z(\dot{r}+\dot{\rho})}{(r+\rho)\sqrt{(r+\rho)^2 - z^2}}
    \end{aligned}
\end{equation}
where $r=\frac{p}{k}$ and $\dot{r} = \sqrt{\frac{\mu}{p}} e \sin f $. Although scenarios in which $\theta$ is outside of the range $(-\frac{\pi}{2},\frac{\pi}{2})$ are of limited interest, it is best to compute $\theta$ using the four-quadrant inverse, $\texttt{atan2(y,r+x)}$.

\subsection*{Spherical Coordinates to Relative Position and Velocity}
To convert the spherical coordinates $(\rho,\theta,\phi)$ and their derivatives $(\dot{\rho},\dot{\theta},\dot{\phi})$ to relative position and velocity vector components, use the transformations
\begin{equation}
    \begin{aligned}
        x &= (r + \rho) \cos \phi \cos \theta - r \\
        y &= (r + \rho) \cos \phi \sin \theta \\
        z &= (r + \rho) \sin \phi \\
        \dot{x} &= (\dot{r} + \dot{\rho}) \cos \phi \cos \theta - (r + \rho)(\dot{\phi} \sin \phi \cos \theta + \dot{\theta} \cos \phi \sin \theta ) - \dot{r}  \\
        \dot{y} &= (\dot{r} + \dot{\rho}) \cos \phi \sin \theta - (r + \rho)(\dot{\phi} \sin \phi \sin \theta - \dot{\theta} \cos \phi \cos \theta ) \\
        \dot{z} &= (\dot{r} + \dot{\rho}) \sin \phi + (r + \rho) \dot{\phi} \cos \phi \\
    \end{aligned}
\end{equation}
where $r=\frac{p}{k}$ and $\dot{r} = \sqrt{\frac{\mu}{p}} e \sin f $.

\subsection*{To Nondimensional Coordinates}

Relative position and velocity vectors are nondimensionalized by normalizng by the orbit radius and changing the independent variable from time to true anomaly. This is efficiently expressed through the transformations
\begin{equation}
    \begin{aligned}
        \boldsymbol{\updelta} \mathbf{\Tilde{r}} &= \frac{1}{r} \boldsymbol{\updelta} \mathbf{r} \\
        \boldsymbol{\updelta} \mathbf{\Tilde{v}} &= -\frac{e}{p} \boldsymbol{\updelta} \mathbf{r} \sin{f} + \frac{1}{k}\sqrt{\frac{p}{\mu}} \boldsymbol{\updelta} \mathbf{v}
    \end{aligned}
\end{equation}
Similar transformations are used to nondimensionalize the spherical coordinates, but only $\rho$ is normalized by the orbit radius. The transformations therefore become
\begin{equation}
    \begin{aligned}
        \Tilde{\rho} &= \frac{\rho}{r} \\
        \Tilde{\rho}' &= -\frac{e}{p} \rho \sin{f} + \frac{\dot{\rho}}{k}\sqrt{\frac{p}{\mu}} \\
        \theta' &= \frac{ \dot{\theta}}{k^2}\sqrt{\frac{p^3}{\mu}} \\
        \phi' &= \frac{ \dot{\phi}}{k^2}\sqrt{\frac{p^3}{\mu}} \\
    \end{aligned}
\end{equation}

\subsection*{From Nondimensional Coordinates}

To convert the nondimensional relative position and velocity vectors to their dimensional forms, the independent variable is converted from true anomaly to time and the vectors are scaled by the orbit radius. This may be accomplished in a single step using
\begin{equation}
    \begin{aligned}
        \boldsymbol{\updelta} \mathbf{r} &= r \boldsymbol{\updelta} \mathbf{\Tilde{r}} \\
        \boldsymbol{\updelta} \mathbf{v} &= \sqrt{\frac{\mu}{p}} \left( e \boldsymbol{\updelta} \mathbf{\Tilde{r}} \sin{f} + k \boldsymbol{\updelta} \mathbf{\Tilde{v}} \right)
    \end{aligned}
\end{equation}
For spherical coordinates, only $\Tilde{\rho}$ is scaled so the transformations become
\begin{equation}
    \begin{aligned}
        \rho &= r \Tilde{\rho} \\
        \dot{\rho} &= \sqrt{\frac{\mu}{p}} \left( e \Tilde{\rho} \sin{f} + k \Tilde{\rho}' \right) \\
        \dot{\theta} &=  \theta' k^2 \sqrt{\frac{\mu}{p^3}} \\
        \dot{\phi} &= \phi' k^2 \sqrt{\frac{\mu}{p^3}} \\
    \end{aligned}
\end{equation}

\newpage
\section*{APPENDIX: Solution Coefficients}

\begin{equation}
    c_{\rho j} = \sum_i \sum_{k\geq i} c_{\rho jik}
\end{equation}

\begin{align*}
    c_{\rho j11} &= \frac{1}{2} K_1^2 \left(1 - 3k_0 \frac{1+2k_0}{1-e^2} \right) \\
    c_{\rho j12} &= -K_1 K_2 \frac{3 + 7k_0}{1-e^2} k_0^2 \sin f_0 \\
    c_{\rho j13} &= K_1 K_3 \frac{2e - (3 + 7k_0) \cos f_0}{1-e^2} k_0^2 \\
    c_{\rho j22} &= K_2^2 \frac{k_0 - 2(1+2k_0)\sin^2 f_0}{1-e^2} k_0^3 \\
    c_{\rho j23} &= -2 K_2 K_3 \frac{1+2k_0}{1-e^2} k_0^3 \sin 2f_0 \\
    c_{\rho j33} &= K_3^2 \frac{e^2 + k_0^2 - 2k_0 (1+2k_0)\cos^2 f_0}{1-e^2} k_0^2 \\
    c_{\rho j55} &= K_5^2 \frac{k_0^2}{1-e^2} \cos 2f_0  \\
    c_{\rho j56} &= -2 K_5 K_6 \frac{k_0^2}{1-e^2} \sin 2f_0 \\
    c_{\rho j66} &= -K_6^2 \frac{k_0^2}{1-e^2} \cos 2f_0
\end{align*}

\newpage
\begin{equation}
    c_{\rho s} = \sum_i \sum_{k\geq i} c_{\rho sik}
\end{equation}

\begin{align*}
    c_{\rho s11} &= \frac{3}{4}K_1^2 \frac{3k_0 + 2k_0^2 + e^2}{k_0(1-e^2)} \sin f_0 \\
    c_{\rho s12} &= K_1 K_2 \frac{6 - 3k_0 + (10 + 7 k_0) \sin^2 f_0}{2(1-e^2)} k_0 \\
    c_{\rho s13} &= K_1 K_3 \frac{e(k_0 - 5) + (10 + 7 k_0)k_0 \cos f_0}{2(1-e^2)} \sin f_0 \\
    c_{\rho s22} &= K_2^2 \frac{9 + k_0 - 2(3 + 2k_0 ) \cos^2 f_0}{2(1-e^2)} k_0^2 \sin f_0 \\
    c_{\rho s23} &= K_2 K_3 \frac{e k_0 (k_0 - 2) + (1 - k_0 + 10 k_0^2 + 2 k_0^3 ) \cos f_0 - 2 k_0^2 (3 + 2 k_0) \cos^3 f_0}{1-e^2} \\
    c_{\rho s33} &= K_3^2 \frac{-2 - e^2(k_0 - 1) + 2k_0 - 5 k_0^2 + k_0^3 + 2k_0^2 (3 + 2k_0)\cos^2 f_0}{2(1-e^2)} \sin f_0 \\
    c_{\rho s55} &= -K_5^2 \frac{\cos 2f_0}{2(1-e^2)}(1+k_0)\sin f_0 \\
    c_{\rho s56} &= K_5 K_6 \frac{\sin 2f_0}{1-e^2}(1+k_0) \sin f_0 \\
    c_{\rho s66} &= K_6^2 \frac{\cos 2f_0}{2(1-e^2)}(1+k_0) \sin f_0
\end{align*}

\newpage

\begin{equation}
    c_{\rho c} = \sum_i \sum_{k\geq i} c_{\rho cik}
\end{equation}

\begin{align*}
    c_{\rho c11} &= \frac{3}{4} K_1^2 \frac{(3+2k_0) \cos f_0 + 3e}{1-e^2} \\
    c_{\rho c12} &= K_1 K_2 \frac{(10 + 7k_0) \cos f_0 + 10e}{2(1-e^2)} k_0 \sin f_0 \\
    c_{\rho c13} &= K_1 K_3 \left(\frac{5}{2} - \frac{10 + 7 k_0}{2(1-e^2)} k_0 \sin^2 f_0 + \frac{15}{2(1-e^2)}k_0^2 \right) \\
    c_{\rho c22} &= -K_2^2 \frac{e^3 + 2(3+2k_0)k_0^2\cos^3 f_0 + 2e(1-3k_0^2) + (1+k_0 -11 k_0^2 + 3k_0^3) \cos f_0}{2(1-e^2)} \\
    c_{\rho c23} &= 2 K_2 K_3 \frac{(1-e^2) - 3 k_0 (1-k_0) + k_0(3+2 k_0) \cos^2 f_0}{1-e^2} k_0 \sin f_0 \\
    c_{\rho c33} &= K_3^2 \frac{e k_0 (4-5k_0) + (-1 + 3 k_0 -7 k_0^2 + 5 k_0^3 ) \cos f_0 + 2(3+2k_0)k_0^2 \cos^3 f_0}{2(1-e^2)} \\
    c_{\rho c55} &= -K_5^2 \frac{e + (1 + k_0)\cos f_0}{2(1-e^2)} \cos 2f_0 \\
    c_{\rho c56} &= K_5 K_6 \frac{e + (1+k_0) \cos f_0}{1-e^2} \sin 2f_0 \\
    c_{\rho c66} &= K_6^2 \frac{e + (1+k_0) \cos f_0)}{2(1-e^2)} \cos 2f_0
\end{align*}

\newpage
\section*{APPENDIX: Correction to Slightly-Eccentric Solution}

In the derivation of the slightly-eccentric solution in spherical coordinates by Butcher et al., the dimensionless first-order equations of relative motion are given in their notation as\cite{ButcherBurnettCurv}
\begin{equation*}
    \begin{aligned}
        \delta r'' - 2 \rho \theta' \delta \theta' - \theta'^2 \delta r &= \frac{2}{\rho^3} \delta r \\
        \delta \theta'' + 2 \delta r' \theta' - \frac{\theta''}{\theta'} \delta \theta' + \theta'' \delta r &= 0 \\
        \delta \phi'' - \frac{\theta''}{\theta'} \delta \phi' + \theta'^2 \delta \phi &= 0
    \end{aligned}
\end{equation*}
The leading-order corrections for the effect of eccentricity are found by substituting the approximations
\begin{equation*}
    \begin{aligned}
        \rho = \frac{r_c}{a_c} &\approx 1 - e \cos M + \frac{1}{2}e^2 (1- \cos 2M) + \cdots \\
        \theta' &\approx 1 + 2 e \cos M + \frac{5}{2}e^2 \cos 2 M + \cdots 
    \end{aligned}
\end{equation*}
into the above equations of motion and linearizing for small $e$. This leads to the system of equations presented in that work,
\begin{equation*}
    \begin{aligned}
        \delta r'' - 2 \delta \theta' - 3 \delta r &= e \cos M \left(10 \delta r + 2 \delta \theta' \right) \\
        \delta \theta'' + 2 \delta r' &= e \sin M \left( 2 \delta r - 2 \delta \theta' \right) - 4 e \delta r' \cos M \\
        \delta \phi'' + \delta \phi &= -4 e \delta \phi \cos M - 2 e \delta \phi' \sin M
    \end{aligned}
\end{equation*}
However, the first-order equations of motion contain a small error. The correct equations are
\begin{equation*}
    \begin{aligned}
        \delta r'' - 2 \rho \theta' \delta \theta' - \theta'^2 \delta r &= \frac{2}{\rho^3} \delta r \\
        \delta \theta'' + \frac{2}{\rho} \delta r' \theta' - \frac{\theta''}{\theta'} \delta \theta' + \frac{\theta''}{\rho} \delta r &= 0 \\
        \delta \phi'' - \frac{\theta''}{\theta'} \delta \phi' + \theta'^2 \delta \phi &= 0
    \end{aligned}
\end{equation*}
The introduction of $\rho$ in the term $\frac{2}{\rho} \delta r' \theta'$ of the $\delta \theta''$ equation contributes to the leading-order effect of eccentricity. After the appropriate substitutions and linearizations, the system of equations becomes 
\begin{equation*}
    \begin{aligned}
        \delta r'' - 2 \delta \theta' - 3 \delta r &= e \cos M \left(10 \delta r + 2 \delta \theta' \right) \\
        \delta \theta'' + 2 \delta r' &= e \sin M \left( 2 \delta r - 2 \delta \theta' \right) - 6 e \delta r' \cos M \\
        \delta \phi'' + \delta \phi &= -4 e \delta \phi \cos M - 2 e \delta \phi' \sin M 
    \end{aligned}
\end{equation*}
Note that the only difference is the coefficient of $\delta r'$ on the right-hand side of the $\delta \theta''$ equation. Along with the second-order terms, this system can be solved by the method described in the paper. Although the change to the equations of motion is small, it has a significant impact on the solution accuracy, and the corrected solution was used for comparison in Figures~\ref{fig:ecc_000101} through \ref{fig:sep_011010_001}. Additional errors are present in the higher-order terms of the $\delta \theta''$ equation in the earlier work, but these are not relevant to the solution used for comparison in this paper.

\end{document}